\documentclass[preprint,12pt]{elsarticle}
\usepackage{amssymb,amsmath,amsfonts,amsthm,latexsym,url,enumerate}
\begin{document}

\begin{frontmatter}

\newtheorem{theorem}{Theorem} %
\newtheorem{lemma}[theorem]{Lemma} %
\newtheorem{claim}[theorem]{Claim} %
\newtheorem{corollary}[theorem]{Corollary} %
\newtheorem{proposition}[theorem]{Proposition}
\newtheorem{question}[theorem]{Question} %
\newtheorem{definition}[theorem]{Definition} %

\theoremstyle{remark} %
\newtheorem*{remark}{Remark} %

\def\C{{\mathbb C}}
\def\F{{\mathbb F}}
\def\N{{\mathbb N}}
\def\Z{{\mathbb Z}}
\def\Q{{\mathbb Q}}
\def\Int{{\operatorname{Int}}}

\def\({\left(}
\def\){\right)}

%%%%%%%%%%%%% Top Matter %%%%%%%%%%%%%

\title{Algebraic Cayley Graphs over Finite Fields\tnoteref{t1}}

\author[lu]{M. Lu}
\ead{mlu@math.tsinghua.edu.cn}

\author[wan]{D. Wan}
\ead{dwan@math.uci.edu}

%    Information for third author
\author[wang]{L.-P. Wang\corref{cor1}}
\ead{wangliping@iie.ac.cn}

\author[zhang]{X.-D. Zhang}
\ead{xiaodong@sjtu.edu.cn}

\cortext[cor1]{Corresponding author}

\address[lu]{Department of Mathematical Sciences, Tsinghua University, Beijing,  China}

\address[wan]{Department of Mathematics, University of California, Irvine, CA 92697-3875, USA}

\address[wang]{Institute of Information Engineering, Chinese Academy of Sciences, Beijing,   China}

\address[zhang]{Department of Mathematics, and MOE-LSC, Shanghai Jiao Tong University, \\
Shanghai,  China}
%\date{}
%\maketitle

\begin{abstract}
A new algebraic Cayley graph is constructed using finite fields. It provides a more flexible source of
expander graphs. Its connectedness, the number of connected components, and diameter bound are studied via
Weil's estimate for character sums. Furthermore, we study the algorithmic problem of computing the number of connected
components and establish a link to the integer factorization problem.
\end{abstract}


\end{frontmatter}

\section{Introduction}
For a subset $S$ of a finite abelian group $\Gamma$, the Cayley graph $Cay(\Gamma, S)$ is
the directed graph with vertex set $\Gamma$, and edge set $\{b_1\rightarrow b_2|b_1-b_2\in S\}$.
Cayley graphs play a central role in the construction of expander graphs. A randomly chosen Cayley graph
$Cay(\Gamma, S)$ often has good properties with non-trivial probability. However, deterministically constructing one such good graph
is often more difficult. Typically one needs to assume additional
structure on the group $\Gamma$ and its subset $S$. By an algebraic Cayley graph, we mean that $\Gamma$ is
the multiplicative group of a finite commutative ring and $S\subset \Gamma$ is a subset with certain algebraic structure such as a box or an interval in some sense.
The box algebraic structure makes it possible to use powerful tools from number theory to prove conditionally (assuming
some sort of Riemann hypothesis) that an algebraic Cayley graph $Cay(\Gamma, S)$ does have the desired properties if the box is suitably large.
In this way, algebraic Cayley graph provides a rich source of expander graphs.

An important such example is given by Chung \cite{chung} who uses the multiplicative group of a finite extension of a finite field and take the subset to be a line
in certain sense. The advantage to work
with a finite field is that the needed estimate can sometimes be proved using the celebrated Weil bound for curves over finite fields.
In this paper, we  introduce a more general construction using the multiplicative group of a finite field and taking the subset  to be those elements represented by certain primary polynomials.

Let $\mathbb{F}_q$ be a finite field of $q$ elements with characteristic
$p$. Let $f(x)$ be an irreducible polynomial of degree $n>1$ over $\mathbb{F}_q$.
Our group $\Gamma$ will be
$$\Gamma_f=(\mathbb{F}_q[x]/(f(x)))^{*}=(\mathbb{F}_q[\alpha])^*= \mathbb{F}_{q^n}^*, \ \ \alpha=\overline{x}.$$
The group $\Gamma_f$ is cyclic of order $q^n-1$.
A polynomial $g(x) \in \mathbb{F}_q[x]$ of degree $d>0$ is called primary if  $g(x)$ is a power of an
irreducible polynomial. For $1\leq d< n$, let $P_d$ be the set of monic primary polynomials of degree $d$ in $\mathbb{F}_q[x]$. Our subset $S$ will be
$$E_d=\{g(\alpha)|g\in P_d\}\subset \Gamma_f.$$
Note that in the case $d=1$, the subset $E_1=\alpha+\mathbb{F}_q$ is a line in the $n$-dimensional $\mathbb{F}_q$-vector space $\mathbb{F}_{q^n}$.

\begin{definition}
Let $G_d(n,q,\alpha)$ be the Cayley graph $Cay(\Gamma_f, E_d)$ with vertex set $\Gamma_f$
and edge set $\{\beta_1\rightarrow \beta_2|\beta_2/\beta_1\in E_d\}$.
\end{definition}
It is clear that $G_d(n,q,\alpha)$ is
a regular directed graph of order $q^n-1$ and its degree is given by
$$|E_d| =|P_d| =\sum_{k|d} \frac{1}{k} \sum_{s|k}\mu(s)q^{\frac{k}{s}} \sim \frac{q^d}{d},$$
where $\mu$ is the M\"obius function.  It should be noted that the graph $G_d(n,q,\alpha)$ depends not just on $d, n, q$ but also
on the choice of $\alpha$ (that is, the choice of the irreducible polynomial $f(x)$ which is used to present the extension field $\mathbb{F}_{q^n}$).
In the case $d=1$, $G_1(n,q,\alpha)$ reduces to Chung's graph in \cite{chung},
which has been studied extensively, see \cite{cohen} \cite{katz}\cite{wan97}.  In this paper, we study the general $d$ case.
Our proof is more direct and uses Weil's bound for character sums.

Our first result is the following theorem.

\begin{theorem}
Assume that $n<q^{d/2}+1$, then the graph $G_d(n,q,\alpha)$
is connected and its diameter $D$ satisfies the bound
\[D\leq 2\frac{n}{d}+1+\frac{4\frac{n}{d}\log{(n-1)}}{d\log{q}-2\log(n-1)}.\]
\end{theorem}

In the case $d=1$, this reduces to the diameter bound in \cite{chung} and \cite{wan97}.
The above theorem gives a sufficient condition for the graph to be connected.
If $n\geq q^{d/2}+1$, the graph $G_d(n,q,\alpha)$ is not always connected as the answer
depends on the choice of $\alpha$ or the irreducible polynomial $f(x)$. More precisely, we have

\begin{theorem}
If $\ell>1$ is a divisor of the integer $(q^n-1)$ such that
$n\geq 2d+ 2(|P_d| +1)\log_{q}\ell $, then there is at least one
$\alpha \in \mathbb{F}_{q^n}$ of degree $n$ such that the number of connected components of the graph
$G_{d}(n,q,\alpha)$ is divisible by $\ell$.
\end{theorem}

If $q>2$, $q^n-1$ has the obvious divisor $(q-1) >1$. We obtain

 \begin{corollary} Assume that  $q>2$ and
$n\geq 2d+ 2(|P_d| +1)$. Then there is at least one
$\alpha \in \mathbb{F}_{q^n}$ of degree $n$ such that the number of connected components of the graph
$G_{d}(n,q,\alpha)$ is divisible by $(q-1)$. In particular, $G_d(n, q, \alpha)$ is not connected for at least one degree
$n$ element $\alpha$.
\end{corollary}

As $|P_d| \sim q^d/d$, the bound $2d + 2(|P_d|+1) \sim 2q^d/d$ is roughly the square of the bound
$q^{d/2}$ in Theorem 2. This shows that the condition in Theorem 2 is not too far from being sharp.
For the remaining interval where
$$ q^{d/2} +1 \leq  n \leq 2d+ 2(|P_d| +1) \sim 2q^d/d,$$
we have no results on the connectedness of the graph $G_d(n, q, \alpha)$. One does know the following crude combinatorial upper bound for the
number $N_d(n, q, \alpha)$ of connected components of the graph $G_d(n, q, \alpha)$:
$$N_d(n, q, \alpha) \leq {{q^n -1}\over { |P_d|+1 \choose {\lceil \frac{n}{d}\rceil -1}}},$$
see Theorem 14 in section $2$.

For a randomly chosen $\alpha$, the graph $G_d(n, q, \alpha)$ is connected with high probability.
For example, this is the case if $\alpha$
is a primitive root of $\mathbb{F}_{q^n}^*$.
Unfortunately, constructing a primitive root (or even an element of high order) is a well known difficult problem
in computational number theory.
In practical application, the difficulty is how to verify quickly that a given $G_d(n, q, \alpha)$ is connected and more generally
how to quickly compute the number of its connected components, using the sparse input size $(n\log q)^{O(1)}$
of the graph $G_d(n, q, \alpha)$. Ideally, we would like to have a deterministic algorithm with running time bounded by a polynomial in $(n\log q)^{O(1)}$,
to compute the number of connected components. In this direction, we have the following conditional
result.

\begin{theorem}
Assume that the factorization of $q^n-1$ is given. Then one can compute the number of
connected components of $G_{d}(n,q,\alpha)$ in time $(n \log q)^{O(1)}$.
\end{theorem}

It would be of great interest to remove the factorization assumption in the above theorem.
The graph $G_d(n, q, \alpha)$ provides a new source of expander graphs, see the last section.
We have

\begin{theorem} Let $\delta$ be a constant with $0<\delta< 1$. Assume that $(n+d-1)\leq q^{d/2}(1 -\delta)$.
Then  the graph $G_d(n, q, \alpha)$ is an expander graph.
\end{theorem}

{\sl Remarks}. In our construction of the Cayley graph $G_d(n, q, \alpha)$, we took the subset $E_d$ to be the set of
all monic primary polynomials of degree $d$. It is also natural to take the subset to the set of all monic irreducible polynomials
of degree $d$ or the set of all monic irreducible polynomials who degree divides $d$. The resulting graph would have
similar quality asymptotically. However, our choice of the subset in this paper makes the proofs simpler and cleaner with
the results slightly better.

%One can extend the results of this paper by replacing our irreducible
%polynomial $f$ by an arbitrary polynomial $f$ of degree $n$. In this case, the group $\Gamma = (\mathbb{F}_q[x]/(f(x)))^{*}$ may not be cyclic. Similar
%results are expected with similar but somewhat more complicated proofs.

\section{The number of  connected components}

Our key technical tool is the following Weil bound for character sums, see Theorem 2.1
in \cite{wan97}.

\begin{lemma} Let $\chi: \Gamma_f \longrightarrow \mathbb{C}^*$ be a non-trivial character.
Then, we have the estimate
$$|\sum_{g\in P_d}\Lambda(g)\chi(g(\alpha))| \leq (n-1)\sqrt{q^d},$$
where $\Lambda(g)$ is the von-Mangold function and it is equal to the degree of the unique prime factor in $g$.
\end{lemma}

\begin{theorem}
If $n<q^{d/2}+1$, then $G_d(n,q,\alpha)$ is connected.
\end{theorem}

\noindent \textbf{Proof.} If the graph $G_d(n,q, \alpha)$ is not connected, then $E_d$ generates a proper subgroup
$H$ of $\Gamma_f$. Let
\[\chi: \Gamma_f \rightarrow \Gamma_f/H \rightarrow \mathbb{C}^*\]
be a non-trivial character of $\Gamma_f$, trivial on $H$.  Then by the Weil bound in Lemma 7,
\[q^d=\sum_{g\in P_d}\Lambda(g)=|\sum_{g\in P_d}\Lambda(g)\chi(g(\alpha))|\leq (n-1)\sqrt{q^d}.\]
It follows that $n\geq q^{d/2}+1$. \hfill $\Box$

The next result shows that the condition $n<q^{d/2}+1$ in the above theorem is not too far from being sharp.

\begin{theorem}
If $\ell>1$ is a divisor of $(q^n-1)$ such that
$n\geq 2d+ 2(|P_d| +1)\log_{q}\ell $, then there is at least one
$\alpha \in \mathbb{F}_{q^n}$ of degree $n$ over $\mathbb{F}_q$ such that the number of connected components of the graph
$G_{d}(n,q,\alpha)$ is divisible by $\ell$.
\end{theorem}

\noindent \textbf{Proof.} Let $\pi_n$ denote the number of monic irreducible polynomials of degree $n$ in
$\mathbb{F}_q[x]$. It is easy to check that
$$|\pi_n -\frac{q^n}{n} | \leq \frac{1}{n}\sum_{k|n, k\leq n/2}q^k  \leq \frac{2}{n} q^{n/2}.$$
The number of degree $n$ elements in $\mathbb{F}_{q^n}$ is $n\pi_n$.
The number of elements in $\mathbb{F}_{q^n}$ which are in a proper subfield of $\mathbb{F}_{q^n}$ is
$$(\sum_{\deg(\alpha) < n} 1) = |n\pi_n -{q^n} | \leq {2} q^{n/2}.$$
Let $H$ be the subgroup generated by $g(\alpha)$ for $g \in P_d$. It is clear that the number of connected
components of the graph $G_d(n,q, \alpha)$ is equal to the index $[\mathbb{F}_{q^n}^* : H]$.

For a divisor $\ell>1$ of $q^n-1$, let $H_{\ell}$ denote the unique subgroup of index $\ell$ in
the cyclic group $\mathbb{F}_{q^n}^*$. The group $H_{\ell}$  consists of $\ell$-th power of elements
in $\mathbb{F}_{q^n}^*$.
Let $I_d$ denote the set of monic irreducible polynomial $g$ in $\mathbb{F}_q[x]$ such that $\deg(g)$
divides $d$. Every element of $P_d$ is an integral power of an element in $I_d$. Furthermore,
$|I_d| =|P_d|$.
If $\alpha$ is a degree $n$ element in $\mathbb{F}_{q^n}$ such that $g(\alpha) \in H_{\ell}$ for all $g \in I_d$,
then $H$ is a subgroup of $H_{\ell}$ and thus the number of connected components of $G_d(n,q, \alpha)$ is
$$[\mathbb{F}_{q^n}^* : H]=[\mathbb{F}_{q^n}^* : H_{\ell}] [H_{\ell}: H] =\ell [H_{\ell}: H]$$
which is divisible by $\ell$. Let
$$N_{\ell} = \{ \alpha \in \mathbb{F}_{q^n}| \deg(\alpha) =n, \ g(\alpha) \in H_{\ell} ~ \forall g\in I_d\}.$$
To prove the theorem, it is enough to prove that $N_{\ell}>0$.
A standard character sum argument shows that
\begin{eqnarray*}
\ell^{|I_d|}N_{\ell} &=& \sum_{\deg(\alpha)=n} \prod_{g\in I_d} \sum_{\chi_g^{\ell}=1} \chi_g(g(\alpha))  \\
&=& \sum_{\chi_g^{\ell}=1, g\in I_d} \sum_{\deg(\alpha)=n} \prod_{g\in I_d} \chi_g(g(\alpha)),
\end{eqnarray*}
where $\chi_g$ denotes a character of $\mathbb{F}_{q^n}^*$.
In the case that $\chi_g =1$ for all $g\in I_d$, the inner sum is the number $n\pi_n$ of degree $n$ elements
in $\mathbb{F}_{q^n}$. In all other $(\ell^{|I_d|}-1)$ cases, there is at least one $g\in I_d$ such that $\chi_g$ is a nontrivial
character. In such a case, the standard Weil character sum bound (see Corollary 2.3 in \cite{wan97}) implies
\begin{eqnarray*}
 |\sum_{\deg(\alpha)=n} \prod_{g\in I_d} \chi_g(g(\alpha))| &=& | \sum_{\alpha \in \mathbb{F}_{q^n}}\prod_{g\in I_d} \chi_g(g(\alpha))
 -\sum_{\deg(\alpha)<n}\prod_{g\in I_d} \chi_g(g(\alpha))| \\
 &\leq &
 ((\sum_{g\in I_d} \deg(g)) -1) q^{n/2}  +  \sum_{\deg(\alpha)<n}1\\
  &\leq &
 (q^d -1) q^{n/2}  +2 q^{n/2} \\
 &=& (q^d+1)q^{n/2}.
 \end{eqnarray*}
Putting these together, we deduce that
\begin{eqnarray*}
\ell^{|I_d|}N_{\ell} &\geq & n\pi_n - (\ell^{|I_d|}-1) (q^d+1)q^{n/2} \\
                             &\geq & q^n - 2q^{n/2} - (\ell^{|I_d|}-1) (q^d+1)q^{n/2} \\
                             &\geq & q^n - \ell^{|I_d|}(q^d+1)q^{n/2} \\
                             & >& q^{\frac{n}{2} +d} (q^{\frac{n}{2}-d} - \ell^{|I_d|+1}).
\end{eqnarray*}
Solving the inequality
$$q^{\frac{n}{2}-d} \geq  \ell^{|I_d|+1},$$
one obtains the condition
$$n\geq 2d+ 2(|I_d| +1)\log_{q}\ell.$$
Since $|I_d| =|P_d|$, the theorem is proved. \hfill $\Box$

In the case $d=1$, we have $|I_1|=|P_1|=q$. This gives the following result.

\begin{corollary}
If $\ell>1$ is a divisor of the integer $q^n-1$ such that $n \geq 2 + 2(q+1)\log_q\ell$,
then there is at least one degree $n$ element $\alpha$ in $\mathbb{F}_{q^n}^*$
such that the graph $G_1(n, q, \alpha)$ is not connected.
\end{corollary}

The above theorem shows that the graph  $G_{d}(n,q,\alpha)$ is not always connected. It depends very much on the
choice of $\alpha$.  An interesting question is to
find a fast algorithm, with running time bounded by a polynomial in $(n\log q)^{O(1)}$,
to compute the number of connected components. In this direction, we have the following conditional
result.

\begin{theorem}
Assume that the factorization of $q^n-1$ is given. Then one can compute the number of
connected components of $G_{d}(n,q,\alpha)$ in time $(n \log q)^{O(1)}$.
\end{theorem}

\textbf{Proof.}
We may assume that $n\geq q^{d/2}+1$, otherwise $G_d(n,q,\alpha)$ is already connected.
Let
\[q^n-1=p^{k_1}_1\cdots p^{k_s}_s, H_i=\{\beta^{p_i}|\beta\in \mathbb{F}^{*}_{q^n}\}.\]
The $H_i$'s are the maximal subgroups of $\mathbb{F}^{*}_{q^n}$. The graph $G_d(n,q,\alpha)$
is dis-connected if and only if the subgroup $H=<g(\alpha)|g\in P_d>$ is contained in
$H_i$ for some $i$. This is true if and only if
\[g(\alpha)^{(q^n-1)/p_i}=1, \forall g \in P_d.\]
The elements of $P_d$ can be listed in time $q^d (n\log q)^{O(1)}$. Note that
$$\max\{ s, (k_1+\cdots +k_s), q^d\} \leq  n^2\log q.$$
It follows that one can check if there is $1\leq i\leq s$ such that $H \subseteq H_i$ in time
\[s q^d(n\log q)^{O(1)}=(n\log q)^{O(1)}.\]
If $H \not\subseteq H_i$ for $1\leq i \leq s$, then $H =\Gamma_f$ and the graph is connected.
Otherwise, we can assume that $H\subseteq H_i$ for some given $i$.

The group  $H_i$ is cyclic of order
$$\frac{q^n-1}{p_i}=p^{k_1}_1\cdots p^{k_{i}-1}_{i}\cdots p^{k_s}_s.$$
Its maximal subgroups are
$H_{ij}=\{\beta^{p_ip_j}|\beta\in \Gamma_f\}
=\{\beta^{p_j}|\beta\in H_i\}$, where $p_ip_j |(q^n-1)$.
Similarly we have $H\subseteq H_{ij}$ for some $j$ if and only if
\[g(\alpha)^{\frac{q^n-1}{p_ip_j}}=1, \forall g \in P_d.\]
Again, we can check if there is $1\leq j \leq s$ such that $H\subseteq H_{ij}$  in time
\[s q^d(n\log q)^{O(1)}=(n\log q)^{O(1)}.\]
Continuing in this fashion, eventually one finds that
$H=H_{i_1i_2\cdots i_u}$, and thus the number of connected components is $[\Gamma_f:H]=p_{i_1}\cdots p_{i_u}$.
The total time needed is bounded by

\hspace*{2cm} $(k_1+\cdots +k_s) q^d(n\log q)^{O(1)}=(n\log q)^{O(1)}$.   \hfill $\Box$
\begin{corollary}
The number of connected components of $G_{d}(n,q,\alpha)$, which is the index $[\Gamma_f:H]$,
can be computed in time $O(q^{n/4})$.
\end{corollary}

\textbf{Proof.} By the well known LLL lattice factorization algorithm,  $q^n-1$
can be factored in time $O(q^{n/4})$.  \hfill  $\Box$

\begin{corollary}
If $n$ is even, the number of connected components of $G_{d}(n,q,\alpha)$ can be computed
in time $O(q^{n/8})$.
\end{corollary}

\textbf{Proof.}  $q^n-1=(q^{n/2}-1)(q^{n/2}+1)$ can be factored in time $O(q^{n/8})$. \hfill  $\Box$

Let $N_d(n, q, \alpha)$ denote the number of connected components of the graph $G_d(n, q, \alpha)$.
An interesting problem is to give a good general upper bound for $N_d(n, q, \alpha)$, which is uniform in $\alpha$.
In this direction, we have
the following simple crude upper bound.
\begin{theorem}
$$N_d(n, q, \alpha) \leq {{q^n -1}\over { |P_d|+1 \choose {\lceil \frac{n}{d}\rceil -1}}}.$$
\end{theorem}

{\sl Proof}. Let $H$ be the subgroup generated by $\{ g(\alpha) | g\in P_d\}$.  Since $\alpha$ has degree $n$,
the unique factorization of polynomials implies that the elements
$$g_1(\alpha) \cdots g_k(\alpha), 0\leq k \leq \lceil \frac{n}{d}\rceil -1, \{g_1, \cdots, g_k\} \subset P_d$$
are distinct elements of $H$. This proves
$$|H| \geq { |P_d|+1 \choose {\lceil \frac{n}{d}\rceil -1}}.$$
It follows that
$$N_d(n, q, \alpha) = [\Gamma_f: H] = \frac{q^n-1}{|H|} \leq  {{q^n -1}\over { |P_d|+1 \choose {\lceil \frac{n}{d}\rceil -1}}} .$$
The theorem is proved. \hfill  $\Box$

\section{The diameter}

The diameter of $G_{d}(n,q,\alpha)$ is the minimal integer $D$ (or $\infty$ if it does not
exist) such that every element in $\Gamma_f$ can  be written as a product of at most
$D$ elements in $E_d$.
%In addition, we define $D^{+}(t,q,d)=\max_{\alpha}D(t,q,\alpha,d)$
%and $D^{-}(t,q,d)=\min_{\alpha}D(t,q,\alpha,d)$.

\begin{theorem}
Assume that $n<q^{d/2}+1$. The diameter $D$ of $G_d(n,q,\alpha)$ satisfies the inequality
\[D\leq 2\frac{n}{d}+1+\frac{4\frac{n}{d}\log{(n-1)}}{d\log{q}-2\log(n-1)}.\]
\end{theorem}

\textbf{Proof.} Let $\widehat{\Gamma_f}$ be the character group of the multiplicative group $\Gamma_f= \mathbb{F}_{q^n}^*$, which is
the set of homomorphisms from $\Gamma_f$ to $\mathbb{C}^*$.
For integer $k>0$ and $\beta\in \Gamma_f$, let $N_{k}(\beta)$ be
the number of solutions of the equation
\[\beta=g_1(\alpha)g_2(\alpha)\cdots g_{k}(\alpha), g_i\in P_d.\]
It is clear that
$$N_{k}(\beta) = \frac{1}{q^{n}-1}\sum_{g_1,\cdots,g_k\in P_d}
\sum_{\chi\in \widehat{\Gamma_f}} \chi(\frac{g_1(\alpha)\cdots g_{k}(\alpha)}{\beta}).$$
To show that the diameter $D$ is bounded by $k$, it is enough to show that $N_{k}(\beta) >0$ for
all $\beta\in \Gamma_f$.  For our purpose, it is simpler to work with the following weighted sum
$$
M_{k}(\beta) = \frac{1}{q^{n}-1}\sum_{g_1,\cdots,g_k\in P_d}
\Lambda(g_1)\cdots \Lambda(g_k)\sum_{\chi\in \widehat{\Gamma_f}} \chi(\frac{g_1(\alpha)\cdots g_{k}(\alpha)}{\beta}).$$
Note that $N_{k}(\beta) >0$ if and only if $M_{k}(\beta) >0$.
Now, separating the trivial character, we obtain
\begin{eqnarray*}
M_{k}(\beta) &=&
\frac{q^{kd}}{q^n-1}+\frac{1}{q^n-1} \sum_{g_1,\cdots,g_k\in P_d}\Lambda(g_1)\cdots \Lambda(g_k) \sum_{\chi\neq 1} \chi(\frac{g_1(\alpha)\cdots g_{k}(\alpha)}{\beta}) \\
&=& \frac{q^{kd}}{q^n-1}+\frac{1}{q^n-1} \sum_{\chi\neq 1} \chi^{-1}(\beta) (\sum_{g \in P_d} \Lambda(g) \chi(g(\alpha)))^k.
\end{eqnarray*}
Applying the Weil bound in Lemma 7, we deduce that
\[|M_{k}(\beta)-\frac{q^{kd}}{q^n-1}|<(n-1)^k \sqrt{q^{dk}}.\]
In order for $M_k(\beta) >0$ for all $\beta$, it suffices to have the inequality
\[q^{kd}\geq q^n (n-1)^k q^{kd/2},\]
that is,
\[ q^{kd -2n} \geq (n-1)^{2k}.\]
This is satisfied if
\[k\geq  \frac{2n}{d- 2\log_q(n-1)} =      2\frac{n}{d}+\frac{4\frac{n}{d}\log(n-1)}{d\log{q}-2\log(n-1)}.\]
The theorem is proved. \hfill $\Box$

For a proper divisor $d$ of $n$,
we now make some comparisons between Chung's graph $G_1(\frac{n}{d},q^d,\beta)$ and our more general
construction $G_d(n,q,\alpha)$, where $\beta$ is a root of an irreducible polynomial of degree $n/d$ in $\mathbb{F}_{q^d}[x]$
and $\alpha$ is a root of an irreducible polynomial of degree $n$ in $\mathbb{F}_{q}[x]$.
It is clear that both graphs have $q^n-1$ vertices.
Assume that $n<q^{d/2}+1$. In this case, both $G_1(\frac{n}{d},q^d,\beta)$ and $G_d(n,q,\alpha)$ are connected, and their diameter bounds
\[
D_1\leq 2\frac{n}{d}+1+\frac{4\frac{n}{d}\log (\frac{n}{d}-1)}{d\log q-2\log (\frac{n}{d}-1)}, \
D_2\leq 2\frac{n}{d}+1+\frac{4\frac{n}{d}\log (n-1)}{d\log q-2\log (n-1)}\]
are comparable.
But $G_1(\frac{n}{d},q^d,\alpha)$ is $q^d-$regular and $G_d(n,q,\alpha)$
is $|P_d|$-regular, where $|P_d|\sim \frac{q^d}{d}<q^d$. Thus, $G_d(n,q,\alpha)$ can be significantly
better than $G_1(\frac{n}{d},q^d,\alpha)$ if $n<q^{d/2}+1$, since $G_d(n,q,\alpha)$ has far fewer edges.

\begin{corollary}
If $q^d>(n-1)^{4\frac{n}{d}+2}$, then $D\leq 2\frac{n}{d}+1$.
\end{corollary}

If $q$ is sufficiently large, it may be possible to improve the above diameter bound to $D \leq \frac{n}{d} +2$.
This is indeed the case for $d=1$, as shown by Katz \cite{katz} and Cohen \cite{cohen}.

A computational question is to ask for a fast algorithm, with running time bounded by $O(n\log q)^{O(1)}$,
to compute the diameter $D$ of the graph $G_d(n, q, \alpha)$. This is expected to be a very difficult problem. Even assuming the
factorization of $q^n-1$, we still do not know a fast algorithm to compute the diameter. We believe that computing the diameter is related to the discrete logarithm
problem and the subset sum problem, both are difficult problems used in cryptography.

\section{Expander graphs}

In this section, we show that our graph  $G_d(n, q, \alpha)$ has good expanding properties.
The adjacency matrix $M=(m_{\beta_1\beta_2})$ is a $(q^n-1) \times (q^n-1)$ matrix, where
the entry $m_{\beta_1\beta_2} = 1$ if $\beta_1\rightarrow \beta_2$ is an edge and it is zero otherwise.  The adjacency operator $M$ acts
on the $(q^n-1)$-dimensional complex vector space $\mathbb{C}^{\Gamma_f}$ of functions on $\Gamma_f$. If $h(x)$ is a complex function on $\Gamma_f$, then
$$M(h)(x) =\sum_{x\rightarrow y} h(y) = \sum_{g\in P_d} h(xg(\alpha)),$$
where $y$ runs over all elements of $\Gamma_f$ such that $x\rightarrow y$ is an edge of $G_d(n,q, \alpha)$.
If $h(x)=\chi(x)$ is a multiplicative character of $\Gamma_f$, then one checks that
$$M(\chi)(x) = \sum_{g\in P_d} \chi(xg(\alpha))= \lambda_d(\chi) \chi(x),$$
where
$$\lambda_d(\chi) = \sum_{g\in P_d} \chi(g(\alpha)).$$
This shows that each character $\chi$ is an eigenvector of the operator $M$. By Artin's lemma, the set of characters on $\Gamma_f$
is $\mathbb{C}$-linearly independent. Since the number of characters is equal to $q^n-1$, it follows that $\mathbb{C}^{\Gamma_f}$ has a
basis consisting of the eigenvectors $\chi$ of $M$, where $\chi$ runs through all characters of $\Gamma_f$.  If $\chi$ is
a character which is trivial on the subgroup generated by $H=<g(\alpha) |g \in P_d>$ of $\Gamma_f$, then the eigenvalue
$$\lambda_d(\chi) = \sum_{ g\in P_d} 1  = |P_d|$$
which is the trivial eigenvalue $\lambda_{triv}=|P_d|$. If $\chi$ is a character which is non-trivial on $H$, its eigenvalue is called a non-trivial eigenvalue which
satisfies the bound
\begin{eqnarray*}
|\lambda_d(\chi)| &= &  |\sum_{g\in P_d} \chi(g(\alpha))| \\
&=& |\frac{1}{d}\sum_{g\in P_d} \Lambda(g)\chi(g(\alpha)) + \sum_{g\in P_d, \Lambda(g)<d}(1 -\frac{\Lambda(g)}{d}) \chi(g(\alpha))| \\
&\leq & \frac{n-1}{d} q^{d/2} +  \sum_{g\in P_d, \Lambda(g)<d} (1 -\frac{\Lambda(g)}{d})     \\
&\leq &  \frac{n-1}{d} q^{d/2} + q^{d/2} \leq \frac{n+d-1}{d}q^{d/2}.\\
\end{eqnarray*}
Since
$$\frac{q^d}{d} =\sum_{g \in P_d} \frac{\Lambda(g)}{d} \leq \sum_{g\in P_d} 1 = |P_d|, $$
we deduce

\begin{theorem} Let $\delta$ be a constant with $0<\delta< 1$. Assume that $(n+d-1)\leq q^{d/2}(1 -\delta)$. Then each non-trivial
eigenvalue $\lambda$  of the adjacency operator $M$ for the graph $G_d(n,q, \alpha)$ satisfies the bound
$$|\lambda| \leq \frac{q^{d}}{d}(1-\delta)  \leq |P_d|(1-\delta) =\lambda_{triv}  (1-\delta).$$
In particular, the graph $G_d(n, q, \alpha)$ is an expander graph.
\end{theorem}

Note that the number of connected components of $G_d(n, q, \alpha)$ is equal to the multiplicity of the trivial eigenvalue $|P_d|$ of the adjacency matrix $M$.
If one uses the matrix $M$ and linear algebra directly to compute the number of connected components, then the running time will be $O(q^n)^{O(1)}$,
which is fully exponential in terms of $n\log q$. This trivial algorithm is far slower than the conditional result in Theorem 5.

Finally, we explain that it is best to view  our graph  $G_d(n, q, \alpha)$ as a weighted graph. For this purpose, let $G_d^*(n,q,\alpha)$ be the weighted graph
with the same vertices and edges as $G_d(n,q,\alpha)$. Given an edge $\beta_1\rightarrow \beta_2$ in $G_d^*(n, q, \alpha)$, we define
the weight of the edge $\beta_1\rightarrow \beta_2$ to be $\Lambda(\beta_2/\beta_1)= \Lambda(g)$, where $\beta_2/\beta_1 = g(\alpha)$
for a unique monic primary polynomial $g \in P_d$. The weighted adjacency matrix $M^*=(m_{\beta_1\beta_2})$ is a $(q^n-1) \times (q^n-1)$ matrix, where
the entry $m_{\beta_1\beta_2} = \Lambda(\beta_2/\beta_1)$ if $\beta_1\rightarrow \beta_2$ is an edge and it is zero otherwise.  The adjacency operator $M^*$ acts
on the $(q^n-1)$-dimensional complex vector space $\mathbb{C}^{\Gamma_f}$ of functions on $\Gamma_f$. If $h(x)$ is a complex function on $\Gamma_f$, then
$$M^*(h)(x) =\sum_{x\rightarrow y} \Lambda(\frac{y}{x}) h(y) = \sum_{g\in P_d} \Lambda(g) h(xg(\alpha)),$$
where $y$ runs over all elements of $\Gamma_f$ such that $x\rightarrow y$ is an edge of $G_d^*(n,q, \alpha)$.
If $h(x)=\chi(x)$ is a multiplicative character of $\Gamma_f$, then one checks that
$$M^*(\chi)(x) = \sum_{g\in P_d} \Lambda(g) \chi(xg(\alpha))= S_d(\chi) \chi(x),$$
where
$$S_d(\chi) = \sum_{g\in P_d} \Lambda(g)\chi(g(\alpha)).$$
This shows that each character $\chi$ is an eigenvector of the operator $M^*$.  If $\chi$ is
a character which is trivial on the subgroup generated by $H=<g(\alpha) |g \in P_d>$ of $\Gamma_f$, then the eigenvalue
$$S_d(\chi) = \sum_{ g\in P_d} \Lambda (g)  = q^d$$
which is the trivial eigenvalue $\lambda_{triv}=q^d$. If $\chi$ is a character which is non-trivial on $H$, its eigenvalue is called a non-trivial eigenvalue which
satisfies the bound
$$|S_d(\chi)| =   |\sum_{g\in P_d} \Lambda(g)\chi(g(\alpha))| \leq (n-1)\sqrt{q^d}. $$
We obtain

\begin{theorem} Let $\delta$ be a constant with $0<\delta< 1$. Assume that $(n-1) \leq q^{d/2}(1 -\delta)$. Then each non-trivial
eigenvalue $\lambda$  of the adjacency operator $M^*$ for the weighted graph $G_d^*(n,q, \alpha)$ satisfies the bound
$$|\lambda| \leq \lambda_{triv}  (1-\delta).$$
In particular, the weighted graph $G_d^*(n, q, \alpha)$ is an expander graph.
\end{theorem}

The condition $(n-1) \leq q^{d/2}(1 -\delta)$ in this weighted theorem is weaker and simpler than the condition $(n+d-1)\leq q^{d/2}(1 -\delta)$ in the previous
un-weighted theorem.

\section*{Acknowledgement}
{\small
 The research is partially supported by
National Natural Science Foundation of China
(Grant No. 10990011, 61170289 and 11271256).}

\end{document}